\documentclass[a4paper,10pt]{elsarticle}
 
\usepackage{amssymb}

\usepackage{graphicx}

\usepackage{epsfig}

\usepackage{subfigure}

\usepackage{color}

\usepackage{amsmath, amssymb, latexsym, amscd}
\usepackage{changepage}
\usepackage{pstricks,pst-plot,pst-node,pst-grad,pst-3d}

\journal{Advances in Water Resources}

\newtheorem{theorem}{Theorem}[section]

\newtheorem{proposition}[theorem]{Proposition}

\newtheorem{example}[theorem]{Example}

\newcommand{\norm}[1]{|\!|#1|\!|}

\def\dis{\displaystyle}

\def\cT{{\mathcal T}}

\def\cT{\mathcal{T}}

\def\begp{\begin{pmatrix}}
\def\ep{\end{pmatrix}}
\def\begm{\begin{matrix}}
\def\edma{\end{matrix}}

\def\dis{\displaystyle}

\def\b0{{\mathbf 0}}

\def\bn{\mathbf{n}}

\def\hx{{\hat x}}

 \newsavebox{\savepar}

\begin{document}

\begin{frontmatter}

\title{An IMPES scheme for a two-phase flow in heterogeneous porous media using a structured grid}

\tnotetext[label1]{}
\author{Gwanghyun Jo}
\author{Do Y. Kwak\corref{cor1}\fnref{label2}}
 
\cortext[cor1]{Corresponding Author. Tel: 82423502720 Fax: 82423502710  kdy@kaist.ac.kr}
\fntext[label2]{This work is supported by NRF, contract  No.2014R1A2A1A11053889}

\address{291 Daehak-ro(373-1 Guseong-dong), Yuseong-gu, Daejeon 305-701, Republic of Korea}

\begin{abstract}
We develop a numerical scheme for a two-phase immiscible flow in heterogeneous porous media using a structured grid finite element method, which have been successfully used for the computation of various physical applications involving elliptic equations \cite{li2003new, li2004immersed, chang2011discontinuous, chou2010optimal, kwak2010analysis}.
The proposed method is based on the implicit pressure-explicit saturation procedure.
To solve the pressure equation, we use
an IFEM based on the Rannacher-Turek \cite{rannacher1992simple} nonconforming space, which is a modification of the work in \cite{kwak2010analysis} where
`broken' $P_1$ nonconforming element of Crouzeix-Raviart \cite{crouzeix1973conforming} was developed.

For the Darcy velocity, we apply the mixed finite volume method studied in \cite{chou2003mixed, kwak2010analysis} on the basis of immersed finite element method (IFEM).
In this way, the Darcy velocity of the flow can be computed  cheaply (locally) after we solve the pressure equation.
The computed Darcy velocity is used to solve the saturation equation explicitly.
Thus the whole procedure can be implemented in an efficient way using a structured grid which is independent of the underlying heterogeneous porous media.
Numerical results show that our method exhibits optimal order convergence rates for the pressure and velocity variables,
and suboptimal rate for saturation.
\end{abstract}

\begin{keyword}
Heterogeneous porous media, Multiphase flow, Immersed finite element method, Mixed finite volume method.
\end{keyword}

\end{frontmatter}

\section{Introduction}

Multiphase flows in porous media arise in various disciplines including petroleum engineering, groundwater  remediation, etc., see \cite{ewing1983th,peaceman1977fundamentals,bear1972dynamics,helmig1997multiphase,bastian1999numerical,ewing1984convergence,russell1983finite} and references therein.

Various discretization methods of the porous media problems have been developed.
These methods include finite difference/volume methods  \cite{douglas1983finite, russell1983finite, michel2003finite, peaceman1977fundamentals}, and
control volume methods \cite{durlofsky1994accuracy, bear1972dynamics, forsyth1991control}.
{Difficulties in finding accurate numerical solutions of the flows in porous media  can arise if   material properties
change abruptly.
For example, the permeability tensor often has a high variation over the space region or is even discontinuous across different materials.
This makes the accurate approximation of the Darcy velocity  hard, especially for the finite difference method.
Since the saturation equation  is of hyperbolic type with a term including  Darcy velocity, the inaccurate computation of the Darcy velocity may lead to a bad approximation of the saturation.}

 For an accurate approximation of Darcy velocity, mixed finite element/volume methods approach \cite{douglas1983time, chavent1991unified, chavent1986mathematical} were developed to solve the pressure equation.
The mixed methods were  combined with implicit pressure and explicit saturation (IMPES) scheme  \cite{sheldon1959one, stone1961analysis} to solve the problem in an efficient way \cite{durlofsky1993triangle,huber1999multiphase}, where the pressure and the saturation equations are solved separately.

On the other hand,  discontinuous Galerkin methods which were developed for elliptic problems  \cite{wheeler1978elliptic,arnold2002unified} were
successfully applied for these problems in \cite{epshteyn2007fully, epshteyn2009analysis,kou2014upwind}.
For an excellent review for various numerical methods for solving equations involving heterogeneous porous media,
see  \cite{miller1998multiphase} and references therein.

 To resolve the discontinuity across the material interface, meshes aligned with the interface are usually used for finite element method (FEM) and their variations.
However, aligned meshes yield unstructured grid and as a result, the data structure becomes more complicated.
Thus, the development of fast, stable numerical methods with high accuracy  is an important issue for the porous media problem.

Recently, there have emerged a new method of using  uniform grids for elliptic P.D.E.s with discontinuous coefficient.
Immersed finite element methods (IFEM) in which one allows the interface to cut through the elements have been introduced by Z. Li,  T. Lin and Y. Lin and their coworkers  \cite{li2003new, li2004immersed}, and its convergence behavior was investigated in  \cite{chou2010optimal}.
IFEM using Crouzeix-Raviart nonconforming bases including applications to mixed finite volume method (MFVM) was developed by Kwak et al. in \cite{kwak2010analysis}
and IFEM for nonhomogeneous jump case were considered in \cite{chang2011discontinuous}.

One of the main features of the IFEMs is that they can use any reasonable  grids independent of interface.
For example, some of the advantages of uniform grids for the interface problems are the following  \cite{chang2011discontinuous}:
\begin{itemize}
\item No grids generation is necessary which can save computation time especially for moving interface.
\item It requires smaller number of degrees of freedom than fitted grids.
\item Simple data structure can render easier development of fast solver such as multigrid methods or alternating direction implicit method.
\item The linear systems are symmetric positive definite when the original problems are symmetric positive definite.
\end{itemize}

In this work, we introduce a numerical scheme using the IFEM for a two-phase immiscible flow through  heterogeneous porous media
having distinct permeability.
Using the concept of the global pressure of Chavent and Jaffr{\'e} \cite{chavent1986mathematical}, the pressure and the saturation equations are solved separately with the IMPES scheme.

To solve the pressure equation, we adopt the MFVM introduced in \cite{chou2003mixed, kwak2010analysis}, where we integrate the momentum and mass equation directly on each
 element. To discretize the equation, we adopt the IFEM for Crouzeix-Raviart nonconforming element developed in \cite{kwak2010analysis}
to the rectangular element.
By choosing appropriate trial space (Rannacher-Turek $Q_1$ space) and test space (the gradient of  trial space),  then by eliminating the Darcy velocity from the system, we can  obtain a discrete equation
for pressure unknowns only.


 After we solve the pressure equation, the Darcy velocity (in Raviart-Thomas \cite{raviart1977mixed} basis) is computed by
 using the local residual of the pressure equation. 
This approach for computing the Darcy velocity is more efficient than the classical mixed finite element method (MFEM) which yields a saddle point formulation \cite{raviart1977mixed}, and is an alternative to the post-processing method of Arnold et al. \cite{arnold1985mixed}.
After pressure equation is solved, saturation equation is solved by the control volume method \cite{ronghua1987generalized, chou2000mixed} together with the upwind scheme.
Numerical examples show that our scheme gives  optimal order convergence for both pressure and velocity for all examples we tested.
For saturation, we obtain the errors at least $O(h^{1.5})$ in $L^2$-norm and $O(h)$ in $H^1$-norm for the case of zero capillary and around $O(h^{1.5})$ in $L^2$-norm and $O(h^{0.5})$ in $H^1$-norm for the nonzero capillary case.

The rest of the paper is organized as follows.
 The governing equations of two phase immiscible flows in heterogeneous porous media are described in section 2.
In section 3, we present a brief review  of IFEM and MFVM for second order elliptic equations.
In section 4, we describe our version of IMPES for two phase flows in heterogeneous porous media.
Numerical examples are given in  section 5,  and conclusion follows in section 6.

\section{Model problem}
Let $\Omega=\Omega^+ \cup \Omega^-$  be a polygonal domain in $\mathbb{R}^2$ separated by a $C^1$-interface $\Gamma$ (See figure 1).
 We assume the subdomains $ \Omega^+$ and $\Omega^-$ are occupied by heterogeneous media  having  discontinuous permeability tensor  $\mathbf{K}$.
The flows of wetting phase and non-wetting phases, i.e. $\alpha \in \{n, w\}$ in $\Omega$ are described by Darcy's law and the continuity equation for each phase.
Governing equations for two-phase flow in porous media, in the absence of gravity, for time $[0,T]$ is given by
\begin{align}
\frac{\partial(\Phi\rho_{\alpha}S_{\alpha})}{\partial t}+\nabla \cdot \{\rho_\alpha  \mathbf{u}_\alpha \}-\rho_\alpha q_\alpha&=0, \quad \alpha =w,n, \label{governing_equation1}
\end{align}
where $\mathbf{u}_\alpha$ are the Darcy velocities given by
\begin{align}
\mathbf{u}_\alpha=-\lambda_\alpha \mathbf{K} \nabla p_\alpha \label{governing_equation2},
\end{align}
 saturations $S_{\alpha}$ and pressures $p_{\alpha}$ satisfy the relations:
\begin{align*}
S_w+S_n&=1 \\
p_n-p_w&=p_c(S_w).
\end{align*}
Here, $\lambda_{\alpha}$ is mobility, $\rho_{\alpha}$ is density {and $q_{\alpha}$ is source term} for the phase $\alpha \in \{w,n\}$.
 {Mobility $\lambda_{\alpha}$ is defined by $\lambda_{\alpha}=k_{r\alpha}/\mu_{\alpha}$, where $k_{r\alpha}$ is relative permeability  and $\mu_{\alpha}$ is the viscosity}.
$p_c$ is the capillary pressure, $\Phi$ is the porosity and $\mathbf{K}$ is the permeability tensor.
The permeability $\mathbf{K}$ is symmetric positive definite and mobilities $\lambda_{\alpha}$ is nonnegative for both phase $\alpha \in \{w,n\}$.
We consider the incompressible fluid  where $\Phi$ and $\rho_{\alpha}$ is constant {and we neglect the effect of gravity}.
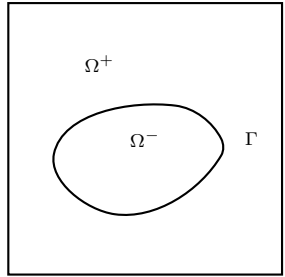
\begin{figure}[t]\label{interface_domain}
\begin{center}
      \psset{unit=2cm}
      \begin{pspicture}(-1,-1)(1,1)
        \pspolygon(0.9,0.9)(-0.9,0.9)(-0.9,-0.9)(0.9,-0.9)
        \psccurve(0.5,0) (0.2,0.22)(-0.6,-0.1)(-0.2,-0.5)(0.5,-0.11)
        \rput(0,0){\scriptsize$\Omega^-$}
        \rput(-0.3,0.5){\scriptsize$\Omega^+$}
        \rput(0.7,0){\scriptsize$\Gamma$}
      \end{pspicture}   \label{fig:doma0}
\caption{A domain $\Omega$ with interface}
\end{center}
\end{figure}

To solve (\ref{governing_equation1}), (\ref{governing_equation2}) we adopt the concept of reformulation by Chavent and Jaffr{\'e} \cite{chavent1986mathematical} where new primary variable global pressure is introduced.
Total velocity $\mathbf{u}$ is defined by sum of $\mathbf{u}_{w}$ and $\mathbf{u}_n$
\begin{align*}
\mathbf{u}&=\mathbf{u}_{n}+\mathbf{u}_w =-\lambda \mathbf{K}(\nabla p_n - f_w \nabla p_c) 
\end{align*}
where total mobility $\lambda$, fractional flow $f_{\alpha}$ are defined by
\begin{align*}
\lambda =& \lambda_w+\lambda_n, \quad f_{\alpha} = \frac{\lambda_{\alpha}}{\lambda} .
\end{align*}
We assume that $\lambda$ is strictly positive.
In \cite{chavent1986mathematical}, global pressure is defined in such a way that:
\begin{align*}
\nabla p=&\nabla p_n - f_w \nabla p_c,
\end{align*}
leading to the relation of total velocity and global pressure:
\begin{align*}
\mathbf{u}&=-\lambda K \nabla p.
\end{align*}
The model problem can be written in an alternative way with unknowns $p$ and $S_w$
\begin{align}
-\nabla \cdot (\lambda \mathbf{K} \nabla p) &= q_w + q_n,   \label{continuous_problem1} \\
\Phi\frac{ S_w}{\partial t}  &=  q_w - \nabla \cdot \left( f_w \mathbf{u} + \lambda_nf_w\mathbf{K}\nabla p_c \right). \label{continuous_problem2}
\end{align}
Initial conditions and boundary conditions are given as follows:
\begin{align*}
\begin{array}{ll}
S_w(\mathbf{x},0) = S_{w0}(\mathbf{x}), \qquad \mathbf{x} \in \Omega & \\
p(\mathbf{x},t)=p_{d}(\mathbf{x},t) \quad \rm{on} \quad \Gamma _{d} &  \mathbf{u} \cdot \mathbf{n} = U(\mathbf{x},t) \quad \rm{on} \quad \Gamma_{n} \\
S_w(\mathbf{x},t)=S_{wd}(\mathbf{x},t) \quad \rm{on} \quad \Gamma _{wd} & \mathbf{u}_w \cdot \mathbf{n} = U_w(\mathbf{x},t) \quad \rm{on} \quad \Gamma_{wn}
\end{array}
\end{align*}

 The goal of this work is to provide numerical methods for the problem (\ref{continuous_problem1}), (\ref{continuous_problem2}).

\section{MFVM and IFEM for second order elliptic equations}
In this section, we briefly explain MFVM  \cite{chou2003mixed}  and IFEM \cite{kwak2010analysis} for elliptic problems having discontinuous  coefficients.
The IFEM was originally designed for $P_1$ nonconforming element of Crouzeix-Raviart, but we modify it to the case of $Q_1$-nonconforming element of Rannacher-Turek \cite{rannacher1992simple}.

Consider the second order elliptic model problem,
\begin{align}\label{elliptic_model1}
\left\{ \begin{array}{rl}
-\rm{div}( \beta  \nabla p ) &=f \quad \rm{in}\quad \Omega \\
p &= 0 \quad \rm{on}\quad \partial \Omega,
\end{array}
\right.
\end{align}
where $f\in L^2(\Omega)$ and $p \in H^1_0 (\Omega)$ and $\beta$ is symmetric and uniformly positive definite matrix with possibly discontinuous entries on the domain $\Omega$.

\subsection{MFVM}
In the computation of porous media equation, it is important to evaluate Darcy velocity accurately, which is a reason for using the MFEM \cite{brezzi1987efficient, brezzi1985two, arnold1985mixed, raviart1977mixed, brezzi1991mixed,douglas1985global}.
Raviart-Thomas element \cite{raviart1977mixed} and Brezzi-Douglas-Marini element \cite{brezzi1985two} are most widely used ones among many MFEMs.
However, one disadvantage of the MFEM formulation is that it leads to a saddle point problem.
One way to avoid the saddle point formulation is the post processing technique of Arnold and Brezzi \cite{arnold1985mixed}.
An alternative, in view of finite volume methods, was suggested  by Chou, Kwak and Kim \cite{chou2003mixed} where the equations in the mixed form are integrated against the gradient of test functions on each element.
We explain it here.
By introducing the velocity $\mathbf{u}=-\beta\nabla p $,
we rewrite (\ref{elliptic_model1}) in a mixed form
\begin{eqnarray}\label{mixed_form1}
\left\{ \begin{array}{rl}
\mathbf{u}+\beta \nabla p &=0, \quad \rm{in} \quad \Omega \\
\rm{div}\mathbf{u} &=f \quad \rm{in} \quad \Omega \\
p &= 0 \quad \rm{on} \quad \partial \Omega.
\end{array}
\right.
\end{eqnarray}
\def\cT{\mathcal T}
For simplicity, assume $\Omega$ is a rectangle and let ${\mathcal T}_h$ be a triangulation by partitioned by rectangles, the case of triangles can be similarly treated.
 
For the basis functions, we use the rotated-$Q_1$ nonconforming finite elements \cite{rannacher1992simple} for the pressure and the lowest-order Raviart-Thomas elements
\cite{raviart1977mixed} for the velocity.
The local space  for pressure on the reference element is given by
\begin{align*}
N_h(\hat{Q})&={\rm span}\{1,\hx,\hat{y},\hat{x}^2-\hat{y}^2\},
\end{align*}
where the degrees of freedoms are the average values on the edges, i.e.,
$\{ \frac{1}{|\hat{e}|} \int_{\hat{e}} \hat{p}_h ds: \hat{e} \hbox{ is and edge of } \cT_h\}.$
The local space for velocity is defined by
\begin{align*}
V_h(\hat{Q})&=\{\hat{v}:\hat{v}=(a+b\hat{x},c+d\hat{y}),a,b,c,d \in \mathbb{R}   \} \\
\end{align*}
and the actually elements are defined  through the Piola mapping. The global spaces $ V_h(\Omega) $ and $N_h(\Omega)$
are defined in an obvious way.

The MFVM scheme for (\ref{mixed_form1}) is the following: Find $(\mathbf{u}_h, p_h) \in V_h(\Omega) \times N_h(\Omega)$, which satisfies
\begin{align}\label{mfvm1}
\left \{ \begin{array}{rl}
\dis \int_Q (\mathbf{u}_h+\beta\nabla p_h) \cdot  \nabla \chi & =0, \qquad \forall \chi \in N_h(Q), \\
\dis \int_Q \rm{div} \mathbf{u}_h & = \dis \int_Q f
\end{array}
\right.
\end{align}
for every element $Q \in \cT_h$.
This  system has an unique solution $(\mathbf{u}_h, p_h) \in V_h(\Omega) \times N_h(\Omega)$ and the following
optimal error estimates holds \cite{chou2003mixed}
\begin{align*}
&\norm{p-p_h}_{L^2(\Omega)}+h|p-p_h|_{H^1(\Omega)} \leq Ch^2\norm{f}_{H^1(\Omega)} \\
&\norm{\mathbf{u}-\mathbf{u}_h}_{L^2(\Omega)}+\norm{{\rm div} \mathbf{u}- {\rm div} \mathbf{u}_h}_{L^2(\Omega)}
\leq Ch(\norm{\mathbf{u}}_{H^1(\Omega)}+\norm{f}_{H^1(\Omega)}).
\end{align*}
Even though the equation (\ref{mfvm1}) is hard to implement, a nice thing about it is that we can transform it to the primal form for $p_h$ only.
Since $\mathbf{u}_h \cdot \bn$ is constant on the edges and $\chi \in N_h(\Omega)$ has common averages on the interior edges, we obtain for each $\chi \in N_h(\Omega)$,
\begin{align*}
\sum_{Q \in \cT_h} \int_Q \mathbf{u}_h \cdot \nabla \chi = \sum_{Q \in \cT_h}  \left[ \int_{\partial Q}(\mathbf{u}_h \cdot n) \chi -\int_Q \rm{div} \mathbf{u}_h \chi \right]=-\int_Q (\overline{f}_h|_Q) \chi,
\end{align*}
where $\overline{f}_h|_Q$ is the local average of $f$ on $Q$.
Considering (\ref{mfvm1}) we obtain
\begin{eqnarray}\label{mfvm_p}
\sum_{Q\in \cT_h}\int_Q \beta \nabla p_h \cdot \nabla \chi = \int_{\Omega}  \overline{f}_h \chi, \quad \chi \in N_h
\end{eqnarray}
which is a symmetric, positive definite system that can be solved easily by conjugate gradient (CG) or preconditioned conjugate gradient (PCG) methods.

 Once $p_h$ is obtained, the velocity $\mathbf{u}_h$ can be computed locally.
Let $Q\in \cT_h$ has edges $e_i$, $i=1,2,3,4$ and $\phi_i \in N_h(Q)$ be the basis function associated with the edge $e_i$.
Then the normal components of velocity are computed locally as   follows\cite{chou2003mixed}:
\begin{align}
|e_i|(\mathbf{u}_h \cdot \bn)|_{e_i}= \int_Q\overline{f}_h \phi_i - \int_Q \beta \nabla p_h \cdot \nabla \phi_i . \label{mfvm_u}
\end{align}

We introduce function spaces, norms, etc, necessary for error analysis.
For any domain $D$, and $m=0,1,\cdots,$ we let $W^m_p(D)$ be the usual Sobolev space with norms(semi-norms)   denoted by
$\norm{\cdot}_{m,p,D}$ ($|\cdot|_{m,p,D}$).
When $p=2$, we write $H^m(D):=W^m_2(D)$ with the norms (semi-norms) $\norm{\cdot}_{m,D}$ ($|\cdot|_{m,D}$).
Let $H^1_0(\Omega)$ be the subspace of $H^1(\Omega)$ with zero trace on the boundary.
 We also define, for $m=1,2$ 
\begin{align*}
\widetilde{H}^m(D) := \{u \in H^{m-1}(D): u|_{D \cap \Omega^s} \in {H}^m(D \cap \Omega^s), s=+,- \}
\end{align*}
with norms,
\begin{align*}
|u|_{m,D}^2 &:= |u|_{m,D\cap \Omega^+}^2 + |u|_{m,D\cap \Omega^-}^2 \\
\norm{u}_{m,D}^2 &:= \norm{u}_{m,D\cap \Omega^+}^2 + \norm{u}_{m,D\cap \Omega^-}^2.
\end{align*}
We let
$H_h(\Omega) :=H^1_0(\Omega) + {N}_h (\Omega)$ equipped with norms $|u|_{1,h}:=(\sum_{Q}|u|^2_{1,Q})^{1/2}$ and $\norm{u}_{1,h}:=(\sum_{Q}\norm{u}^2_{1,Q})^{1/2}$.

\subsection{IFEM for problems with discontinuous coefficients}
In this subsection, we briefly review the IFEM for problem with discontinuous coefficients.
This method allows us to use the grid that does not align with the interface.
 We restrict our attention to the case of scalar $\beta$. Tensor coefficient can be treated similarly.

 Assume
\begin{eqnarray*}
\beta=\left \{ \begin{array}{ll}
\beta^+, & \mathbf{x} \in \Omega^+, \\
\beta^-, & \mathbf{x} \in \Omega^-,
\end{array}
\right.
\end{eqnarray*}
with jump conditions
\begin{eqnarray*}
[p]=0, \left[\beta\frac{\partial p}{\partial n} \right]=0 \quad \rm{across} \quad \Gamma.
\end{eqnarray*}

The IFEMs have been suggested and proven to be efficient for elliptic problems, see the work of Z. Li, T. Lin and Y. Lin, etc., in \cite{li2003new, li2004immersed} and S. Chou, D. Y. Kwak and K. T. Wee in \cite{chou2010optimal, kwak2010analysis}.
In particular, the edge based Crouzeix-Raviart $P_1$-nonconforming IFEM combined with Raviart-Thomas MFVM  \cite{kwak2010analysis} turns out to be efficient in implementing mixed formulation through a primary
form.

We develop the $Q_1$-nonconforming version of IFEM scheme. 
To do so, we modify Rannacher-Turek element \cite{rannacher1992simple} $Q_1$-nonconforming basis function
so that the flux is weakly continuous along the local interface in a manner similar to \cite{kwak2010analysis}. 
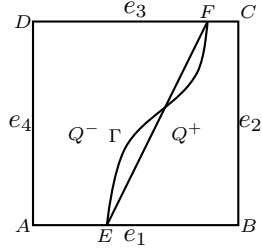
\begin{figure}[ht]
  \begin{center}
    \psset{unit=2.7cm}
    \begin{pspicture}(0,0)(1,1)
      \psset{linecolor=black} \pspolygon(0,0)(1,0)(1,1)(0,1) \psline(0.36,0)(0.85,1)
      \pscurve(0.36,0)(0.45,0.38)(0.8,0.75)(0.85,1)
      \rput(-0.05,0){\scriptsize$A$}
      \rput(1.05,0){\scriptsize$B$}
      \rput(1.05,1.05){\scriptsize$C$}
      \rput(-0.05,1){\scriptsize$D$}
      \rput(0.5,-0.06){$e_1$}
      \rput(1.06,0.5){$e_2$}
      \rput(0.5,1.06){$e_3$}
      \rput(-0.06,0.5){$e_4$}
      \rput(0.85,1.05){\scriptsize$F$}
      \rput(0.25,0.45){\scriptsize$Q^-$}
      \rput(0.75,0.45){\scriptsize$Q^+$}
      \rput(0.35,-0.05){\scriptsize$E$}
      \rput(0.40,0.44){\scriptsize$\Gamma$}
\pnode(-.3,0.6){a}
\pnode(0.12,0.5){b}
    \end{pspicture}
    \caption{A typical interface rectangular} \label{fig:interel}
\end{center}
\end{figure}
To do so, we begin need to explain the interface/element relation:
We call an element $Q \in \cT_h$ is an interface element if the interface $\Gamma$ passes though the interior of $Q$
(see  Figure \ref{fig:interel}).
Suppose interface cut through the edges of $Q$ at point $E$ and $F$ and divide $Q$ by two section $Q^+$ and $Q^-$.
Let $G$ be the midpoint of $\overline{EF}$
Let $e_i$, $i=1,2,3,4,$ be the edges of $Q$.
For $\phi \in H^1(T)$, let $\overline{\phi}_{e_i}$ denote the average of $\phi$ along $e_i$, i.e.,
\begin{eqnarray*}
\overline{\phi}_{e_i} :=\frac{1}{|e_i|}\int_{e_i}\phi ds.
\end{eqnarray*}
We construct a basis function of the form
\begin{eqnarray}\label{im_basis}
\phi =\left \{ \begin{array}{ll}
\phi^+ =a_0+b_0x+c_0y+d_0(x^2-y^2), & (x,y) \in Q^+ \\
\phi^- =a_1+b_1x+c_1y+d_1(x^2-y^2), & (x,y) \in Q^-
\end{array}
\right.
\end{eqnarray}
satisfying
\begin{align}
\overline{\phi}_{e_i}&=V_i, \quad i=1,2,3,4,      \label{im_basis_cond1} \\
\phi^+(E)&=\phi^-(E), \quad \phi^+(F)=\phi^-(F),  \label{im_basis_cond2} \\
d_0&=d_1,\label{im_basis_cond3} \\
\quad\beta^+\frac{\partial\phi^+}{\partial n_{\overline{EF}}}(G)&=\beta^-\frac{\partial \phi^-}{\partial n _{\overline{EF}}}(G),\label{im_basis_cond4}
\end{align}
where $V_i$, $i=1,2,3,4$ are given values, $n_{\overline{EF}}$ is the unit normal vector on the line segment $\overline{EF}$.
It is easy to see the conditions (\ref{im_basis_cond1}), (\ref{im_basis_cond2}), (\ref{im_basis_cond3}) and (\ref{im_basis_cond4}) determine $\phi$ uniquely.
Let $\widehat{N}_h(Q)$ denote the (local) four-dimensional space spanned by these shape functions.

We define the immersed finite element space $\widehat{N}_h(\Omega)$ as the collection of the functions such that
\begin{eqnarray*}
\left \{ \begin{array}{ll}
\phi|_Q \in N_h(Q) & \hbox{if $Q$ is a noninterface element}, \\
\phi|_Q \in \widehat{N}_h(Q) & \hbox{if $Q$ is an interface element}, \\
\int_e\phi|_{Q_1}ds=\int_e\phi|_{Q_2}ds & \hbox{if $Q_1$, $Q_2$ are adjacent elements and $e$ is a common edge}, \\
\int_e\phi ds=0 &\hbox{if $e$ is part of the boundary $\partial \Omega$ }. \\
\end{array}
\right.
\end{eqnarray*}

The
IFEM scheme for the equation (\ref{elliptic_model1}) is given as follows: Find $p_h \in \widehat{N}_h(\Omega)$ such that
\begin{eqnarray}\label{ifem_bilear_form}
a_h(p_h,\phi)=(f,\phi), \quad \forall \phi \in \widehat{N}_h (\Omega),
\end{eqnarray}
where
\begin{align}\label{ifem_bilear_form2}
a_h(q,\phi) &= \sum _{Q \in \cT_h} \int_Q  \beta \nabla q \cdot \nabla \phi dx, \quad \forall q, \phi \in  \widehat{H}_h(\Omega),\\
 \widehat{H}_h(\Omega) &:=H^1_0(\Omega) + \widehat{N}_h (\Omega).
\end{align} Here $ \widehat{H}_h(\Omega) $ is equipped with the norm $\norm{\cdot}_{1,h} .$
A slight modification of the proof in \cite{kwak2010analysis} gives the following result.
\begin{proposition}
Suppose $p \in \widetilde{H}^2(\Omega)$ and $p_h \in \widehat{N}_h(\Omega)$ are solutions to (\ref{mixed_form1}) and (\ref{ifem_bilear_form}).
Then the following error estimates holds:
\begin{eqnarray}\label{proposition1}
& &\norm{p-p_h}_{L^2(\Omega)}+h|p-p_h|_{1,h} \leq Ch^2\norm{p}_{\widetilde{H}^2(\Omega)}.
\end{eqnarray}
\end{proposition}
\subsection{MFVM using IFEM}
We apply our IFEM to MFVM introduced in Section 3.1 to compute Darcy velocity.
The method is similar to MFVM introduced in section 3.1., but the $Q_1$-nonconforming finite element space
 ${N}_h(\Omega)$ is replaced by broken $Q_1$-nonconforming space $\widehat{N}_h(\Omega)$ in (\ref{mfvm1}) and we obtain
the following error estimate:
\begin{proposition}
\begin{align*}
&  \norm{p-p_h}_{L^2(\Omega)}+h|p-p_h|_{1,h} \leq Ch^2\norm{p}_{\widetilde{H}^2(\Omega)} \\
& \norm{\mathbf{u}-\mathbf{u}_h}_{L^2(\Omega)}+\norm{ {\rm div} \mathbf{u}- {\rm div} \mathbf{u}_h}_{L^2(\Omega)}
\leq Ch \{ \norm{\mathbf{u}}_{H^1(\Omega)}+\norm{p}_{\widetilde{H}^2(\Omega)}+\norm{f}_{\widetilde{H}^1(\Omega)} \}
\end{align*}
provided $\mathbf{u} \in \mathbf{H}^1(\Omega)$ and $f \in \widetilde{H}^1(\Omega)$.
\end{proposition}

\section{Application to IMPES scheme }
Suppose the permeability tensor given in the form $\mathbf{K}=K\mathbf{I}$, where
\begin{eqnarray*}
K=\left \{ \begin{array}{ll}
K^+, & \mathbf{x} \mbox{ in } \Omega^+, \\
K^-, & \mathbf{x} \mbox{ in } \Omega^-.
\end{array}
\right.
\end{eqnarray*}

We will briefly review the original IMPES introduced by Sheldon (1959) and Stonde and Gardner (1961) \cite{sheldon1959one, stone1961analysis}.
Consider the problem (\ref{continuous_problem1}),(\ref{continuous_problem2}) on the time interval $[0,T]$.
Let $\triangle t >0$ be a timestep and divide $T$ by $N$ time steps.
Denote $t^l=l\triangle t$ for $l=0, 1, 2, ..., N$.
Equations (\ref{continuous_problem1})-(\ref{continuous_problem2}) are solved sequentially:
Given the saturation at time level $t^l$, the pressure $p$ at time $t^{l+1}$ is computed from (\ref{continuous_problem1}) which is a linear equation in $p$. 
Then the saturation at time level $t^{l+1}$ is computed from (\ref{continuous_problem2}) with velocity fields $\mathbf{u}$
obtained at time level $t^{l+1}$.
\vskip 0.4cm
\textbf{IMPES Algorithm}.
\\ For $l=0,\cdots, N-1,$ repeat the following process.
\begin{enumerate}
  \item Solve for $p^{l+1}$ and $\mathbf{u}^{l+1}$ with given $S_w^l$
\begin{eqnarray}\label{impes1}
\left\{
  \begin{array}{rl}
    \nabla \cdot \mathbf{u}^{l+1}&= q_w^{l+1} + q_n^{l+1}   \\
    \mathbf{u}^{l+1}&=-\lambda^l \mathbf{K}\nabla p^{l+1},
      \end{array}
\right.
\end{eqnarray}
where $\lambda^l=\lambda(S_w^l)$.
  \item Compute $S_w^{l+1}$
  \begin{eqnarray}\label{impes2}
\left\{
  \begin{array}{rl}
    \mathbf{u}_w^{l+1}&= f_w^l \mathbf{u}^{l+1} + \lambda_n^l f_w^l\mathbf{K}\nabla p_c(S_w^{l})    \\
    \Phi(S_w^{l+1}-S_w^l)&=\Delta t( q_w^{l+1} - \nabla \cdot \mathbf{u}_w^{l+1} ),
  \end{array}
\right.
\end{eqnarray}
where $f^l_w=f_w(S_w^l)$.
\item $l=l+1$
\end{enumerate}

We will call (\ref{impes1}) the pressure equation, and (\ref{impes2}) the saturation equation.
In our version of IMPES,  we solve the pressure equation by MFVM based on IFEM introduced in Section 3.3 and then solve the saturation equation by the control volume method.

\subsection{Pressure equation }
Since we use saturation $S_w$  at time level $t^l$, $\lambda^l$ is considered as a coefficient of the pressure equation.

Since global pressure and total velocity are continuous along the interface $\Gamma$, we impose the condition
\begin{eqnarray}\label{immcondition}
[p^{l+1}]=0, \left[ K \frac{\partial p^{l+1}}{\partial n} \right]=0 \quad \rm{across} \quad \Gamma.
\end{eqnarray}
This can be effectively solved by the IFEM  described in section 3.2: Find $p_h^{l+1} \in \widehat{N}_h(\Omega)$ satisfying
\begin{eqnarray*}
a_h(p_h^{l+1},\phi)=(q_w^{l+1}+q_n^{l+1},\phi), \quad \forall \phi \in \widehat{N}_h (\Omega),
\end{eqnarray*}
where is the bilinear form $a_h(\cdot,\cdot)$ is same as (\ref{ifem_bilear_form2}) except that the coefficient $\beta$ is replaced by $\lambda^l \mathbf{K}$.

Once the global pressure $p^{l+1}_h$ is solved by IFEM, the total velocity $\mathbf{u}_h^{l+1}$ is computed by local recovery technique of MFVM in Section 3.3.

\subsection{Saturation equation} First we note that the $Q_1$-nonconforming element can have
negative values even though the degrees of freedom are nonnegative.
Since the saturation variables $S_{\alpha}$ cannot have
negative values, it is inappropriate to use the  the $Q_1$-nonconforming element. Hence
we use the $Q_1$-conforming finite element for the approximation of the saturation.

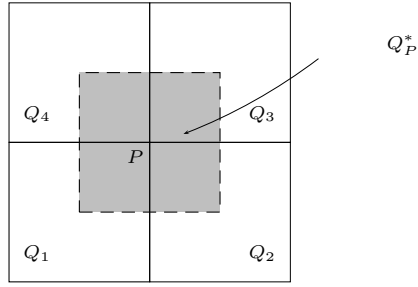
\begin{figure}[ht]
  \begin{center}

    \psset{unit=3.7cm}
    \begin{pspicture}(0,0)(1,1)
      \pnode(0.25, 0.25){a}  \pnode(0.75,0.25){b}  \pnode(0.75, 0.75){c} \pnode(0.25, 0.75){d}
      \psline[linewidth=0.5\pslinewidth, linestyle=dashed, linecolor=black,fillstyle=solid,fillcolor=lightgray](0.25, 0.25)(0.75,0.25)(0.75, 0.75)(0.25, 0.75)(0.25, 0.25)

      \psset{linecolor=black} \psset{linewidth=0.5\pslinewidth} \pspolygon(0,0)(0.5,0)(0.5,0.5)(0,0.5)
      \psset{linecolor=black} \psset{linewidth=0.5\pslinewidth} \pspolygon(0.5,0)(1,0)(1,0.5)(0.5,0.5)
      \psset{linecolor=black} \psset{linewidth=0.5\pslinewidth} \pspolygon(0,0.5)(0.5,0.5)(0.5,1)(0,1)
      \psset{linecolor=black} \psset{linewidth=0.5\pslinewidth} \pspolygon(0.5,0.5)(1,0.5)(1,1)(0.5,1)
      \rput(0.45,0.45){\scriptsize $P$}
      \rput(0.1,0.1){\scriptsize $Q_1$}
      \rput(0.9,0.1){\scriptsize $Q_2$}
      \rput(0.9,0.6){\scriptsize $Q_3$}
      \rput(0.1,0.6){\scriptsize $Q_4$}

      \pnode(1.1,0.8){a} \pnode(0.62, 0.53){b}
      \ncarc[linewidth=0.5\pslinewidth]{->}{a}{b}
      \rput(1.4,.85){\scriptsize $Q_P^*$}
    \end{pspicture}
    \caption{A typical dual volume $Q_P^*$ associated with node $P$.} \label{fig:dualvolume2}
\end{center}
\end{figure}

Let $U_h(\Omega)$ denote $Q_1$-conforming finite element space defined over the triangulation given above.
We apply the control volume method to (\ref{impes2}) together with the upwinding scheme as used in \cite{durlofsky1994accuracy, bear1972dynamics,forsyth1991control}, which
we briefly explain as follows:
Fix an interior vertex $P$ and consider the four elements  $Q_1$, $Q_2$, $Q_3$ and $Q_4$ sharing $P$ as a common vertex. Connecting the centers of the four
elements $Q_i, i=1,\cdots,4$, we obtain a dual volume, which we shall call  $Q_P^*$ (Figure \ref{fig:dualvolume2}).
For a boundary vertex, an obvious modification is necessary.

The  collection of dual volumes $Q_P^*$ is denoted by $\mathcal D_h$.

A test function space $W_h(\Omega)$ using $D_h$ by
\begin{align*}
W_h(\Omega)=\{ \psi \in L^2(\Omega)| \psi \hbox{ is piecewise constant on each } Q_P^* \in \mathcal D_h \}.
\end{align*}

\begin{figure}[ht]
  \begin{center}
    \psset{unit=3.7cm}
    \begin{pspicture}(0,0)(1,1)
      \pnode(0.25, 0.25){a}  \pnode(0.75,0.25){b}  \pnode(0.75, 0.75){c} \pnode(0.25, 0.75){d}
      \psline[linewidth=0.5\pslinewidth, linestyle=dashed, linecolor=black,fillstyle=solid,fillcolor=lightgray](0, 0)(0.5,0)(0.5, 0.5)(0, 0.5)(0, 0)

      \psset{linecolor=black} \psset{linewidth=0.5\pslinewidth} \pspolygon(0,0)(1,0)(1,1)(0,1)
      \rput(-0.06,-.06){\scriptsize$P_i$}
      \rput(1.06,-.06){\scriptsize$P_{j_1^k}$}
      \rput(-0.06,0.94){\scriptsize$P_{j_2^k}$}
      \rput(0.55,0.5){\scriptsize$C_k$}
      \rput(0.55,-.06){\scriptsize$M_{ij}$}
      \rput(0.45,0.15){\scriptsize$\gamma_{ij}^k$}
      \rput(0.6,0.3){\scriptsize$\mathbf{n}_{ij}^k$}
      \rput(0.45,0.65){\scriptsize$Q_k$}
      \pnode(0.5,0.25){a} \pnode(0.7,0.25){b}
      \ncline[linewidth=0.5\pslinewidth]{->}{a}{b}
    \end{pspicture}
    \caption{Illustration of $\gamma_{ij}^k$ in $Q_k$ for $j=j_1^k$.} \label{fig:dualvolume3}
\end{center}
\end{figure}
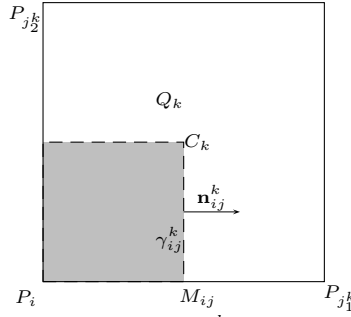

To compute $S_w^{l+1}$ in (\ref{impes2}), we first consider the Petrov-Galerkin methods with test function space $W_h(\Omega)$:
Find $S^{l+1}_w \in U_h$ such that satisfies
\begin{align}\label{control_volume1}
\int_{\Omega} \Phi S^{l+1}_w \psi_i =  \int_{\Omega}(\Phi S^{l}_w+\triangle t q_w^{l+1} )\psi_i  - \int_{\Omega} \triangle t  \nabla \cdot \mathbf{u}_w^{l+1} \psi_i
\end{align}
for all $\psi_i \in W_h(\Omega)$. Fix a vertex $P_i$ shared by four neighboring elements $Q_k, k=1,\cdots,4$   and
applying the divergence theorem in (\ref{control_volume1}) on  $Q_{P_i}^*$  leads to the following form
\begin{align}
\int_{Q_{P_i}^*} S^{l+1}_w &=\int_{Q^*_{P_i}}(S^{l}_w+ \frac{\triangle t}{\Phi} q_w^{l+1} )- \frac{\triangle t}{\Phi} \int_{\partial Q^*_{P_i}}  \mathbf{u}_w^{l+1} \cdot \mathbf{n} ds. \label{control_volume2}
\end{align}
Now we approximate the second term of (\ref{control_volume2}) by the upwinding concept.
Suppose $P_j$ is an adjacent node of $P_i$ belonging to the element $Q_k$ with barycenter $C_k$.
We denote by $M_{ij}$ the mid point of $P_i$ and $P_j$.
Let us denote by $\gamma_{ij}^k$ the segment $\overline{M_{ij}C_k}$, whose outward normal vector is $\mathbf{n}_{ij}^k$.
Then,
\begin{align*}
\int_{\partial Q^*_{P_i} } \mathbf{u}_w^{l+1} \cdot \mathbf{n} ds
=\sum_{k=1}^4 \int_{\partial Q^*_{P_i} \cap Q_k} \mathbf{u}_w^{l+1} \cdot \mathbf{n} ds
=\sum_{k=1}^4 \sum_{j=j_1^k,j_2^k}\int_{\gamma_{ij}^k}  \mathbf{u}_w^{l+1} \cdot \mathbf{n}_{ij}^k ds,
\end{align*}
where $P_{j_1^k}$ and $P_{j_2^k}$ are two  adjacent nodes of $P_i$, belonging to the element $Q_k$ (see Figure \ref{fig:dualvolume3}).
Noting that 
$$\int_{\gamma_{ij}^k}  \mathbf{u}_w^{l+1} \cdot \mathbf{n}_{ij}^k ds= \int_{\gamma_{ij}^k} f_w(S_w)(\mathbf{u}^{l+1}+ {K\lambda_{n}^l \nabla p_c(S_w^l)}          )\cdot  \mathbf{n}_{ij}^k ds,$$
we replace it   by
\begin{align*}
\int_{\gamma_{ij}^k} f_w(S_w^*)(\mathbf{u}^{l+1}+ {K\lambda_{n}^l \nabla p_c(S_w)^l}) \cdot \mathbf{n}_{ij}^k ds,
\end{align*}
where
\begin{align*}
S_w^*=\left\{
\begin{array}{rl}
S_w(P_i), & \hbox{if} \quad (\mathbf{u}^{l+1}+ {K\lambda_{n}^l \nabla p_c(S_w)^l}) \cdot \mathbf{n}_{ij}^k \geq 0 \\
S_w(P_j), & \hbox{if} \quad (\mathbf{u}^{l+1}+ {K\lambda_{n}^l \nabla p_c(S_w)^l}) \cdot \mathbf{n}_{ij}^k < 0.
\end{array}
\right.
\end{align*}

The final scheme for the saturation equation is as follows: Find $S_w^{l+1} \in U_h(\Omega)$ such that:
\begin{align*}
\quad \int_{Q_{P_i}^*} S^{l+1}_w &=\int_{Q^*_{P_i}}(S^{l}_w+ \frac{\triangle t}{\Phi} q_w^{l+1} ) - \frac{\triangle t}{\Phi}  \sum_{k}\sum_{j=j_1^k,j_2^k}\int_{\gamma_{ij}^k}  f_w(S_w^*)(\mathbf{u}^{l+1}+ {K\lambda_{n}^l \nabla p_c(S_w)^l}) \cdot \mathbf{n}_{ij}^k ds,
\end{align*}
for all $P_i$.

\section{Numerical Results}
We present some numerical simulations of two phase immiscible flows in porous media. For the first {three} examples,  the exact solutions
are known and the domain $\Omega$ is $[0, \pi/2]^2$. The {fourth} example is the famous five spot flooding problem with
 $\Omega = (0,300 [ m ]  ) \times (0,300 [ m ] )$.

We  divide the domain by $n \times n$ ($n=2^j (j=2,3,\cdots)$)  rectangles whose sides are parallel to the coordinate  axes.
The relative {permeability} {$k_{r\alpha} $} is given by Brooks-Corey model with index $\lambda=2$ \cite{brooks1966properties}:
\begin{eqnarray*}
\left\{
\begin{array}{l}
k_{rw}(S_w)=S_w^4  \\
k_{rn}(S_w)=(1-S_w)^2(1-S_w^2). \\
\end{array}
\right.
\end{eqnarray*}
In first two examples, we assume zero capillary effects.
In example 5.3, we consider the case of capillary pressure given by Brooks-Corey model with index $\lambda=2$:
\begin{align*}
p_c(S_w)=p_dS_w^{-\frac{1}{2}},
\end{align*}
where $p_d$ is an entry pressure.

We use the  level set function
$$L(x,y): \Omega \rightarrow \mathbb{R}$$
to describe the  interface between the heterogeneities of the domain.
We assume that 
\begin{eqnarray*}
L(x,y)=\left\{
\begin{array}{l}
<0, \quad (x,y) \quad \hbox{in} \quad \Omega^- \\
\quad 0, \quad (x,y) \quad \hbox{on} \quad \Gamma \\
>0, \quad (x,y) \quad \hbox{in} \quad \Omega^+
\end{array}
\right.
\end{eqnarray*}
and the permeability  $K$ is given as follows:
\begin{eqnarray*}
K(x,y)=\left\{
\begin{array}{l}
K^{+}, \quad (x,y) \quad \hbox{in} \quad \Omega^+ \\
K^{-}, \quad (x,y) \quad \hbox{in} \quad \Omega^-.
\end{array}
\right.
\end{eqnarray*}
Our examples are constructed so that fluxes are continuous along the interface $\Gamma$, i.e. $ [\mathbf{u} \cdot \mathbf{n}] = [\mathbf{u}_w \cdot \mathbf{n}]=0$ across $\Gamma$.

\begin{example}\label{example1}

We set the level set function $L(x,y)$ as
\begin{eqnarray*}
L(x,y)=(x+y-1)(x+y-3).
\end{eqnarray*}
Permeability   is given by  $K^+=1$ and $K^-=0.001$.
The exact solutions are  given as follows:  
\begin{align*}
& p(x,y,t)  \\=& \left\{
  \begin{array}{ll}
    (l_{xy}-1)(l_{xy}-2)(2-t)K^{-}+100, & \hbox{$0 \leq l_{xy}<1$} \\
    -(l_{xy}-1)\{1/3(l_{xy}-1)^2-(l_{xy}-1)+1\}(2-t)K^{+}+100, & \hbox{$1 \leq l_{xy} \leq 3$,}  \\
    -\{(l_{xy}-3)(l_{xy}-2)+\frac{2}{3}\}(2-t)K^{-}+100, & \hbox{$3\leq l_{xy} \leq \pi$,}
  \end{array}
\right. \\
\mbox{ and } \\
& S_w(x,y,t) = \left(1-1/8l_{xy}-1/20l_{xy}^2 \right)(0.5 + 0.5t),
\end{align*}
where $l_{xy}=x+y$ and the source terms $q_w$ and $q_n$ are chosen to satisfy   (\ref{continuous_problem1}) and (\ref{continuous_problem2}).

 We set the porosity $\Phi=1$ and the viscosities $\mu_w=\mu_n=1$.  {Dirichlet boundary condition on all of $\partial\Omega$ is imposed.}

Table \ref{Error_Table1} shows the $L^2$, $H^1$  errors of  saturation, pressure and Darcy velocity when  $T=1$.
The time step $\triangle t$ is chosen as  $ \triangle t=16/n^2 =(64/\pi^2) h^2$.
%

\begin{table} \caption{MFVM based on broken $Q_1$-nonconforming IFEM methods for {Example} \ref{example1}}
\begin{scriptsize}
\begin{center}
\begin{tabular}{|c|l c|l c| l c|}
  \hline
  Element  & $\norm{S_w-S_{w,h}}_{L_2(\Omega)}$  & order & $\norm{p-p_{h}}_{L_2(\Omega)}$ & order & $\norm{\mathbf{u}-\mathbf{u}_{h}}_{L_2(\Omega)}$  & order \\
  \hline
  $8^2$  & $1.198 \times 10^{-3}$ &       & $6.068 \times 10^{-2}$ &       & $3.767 \times 10^{-4}$ &       \\
  $16^2$ & $2.994 \times 10^{-4}$ & 2.000 & $1.623 \times 10^{-2}$ & 1.903 & $1.269 \times 10^{-4}$ & 1.569 \\
  $32^2$ & $7.252 \times 10^{-5}$ & 2.046 & $4.012 \times 10^{-3}$ & 2.016 & $4.534 \times 10^{-5}$ & 1.485 \\
  $64^2$ & $1.673 \times 10^{-5}$ & 2.116 & $1.013 \times 10^{-3}$ & 1.987 & $1.937 \times 10^{-5}$ & 1.227 \\
  $128^2$& $3.579 \times 10^{-6}$ & 2.225 & $2.502 \times 10^{-4}$ & 2.017 & $9.176 \times 10^{-6}$ & 1.081 \\
  $256^2$& $8.465 \times 10^{-7}$ & 2.080 & $6.438 \times 10^{-5}$ & 1.958 & $4.521 \times 10^{-6}$ & 1.021  \\  \hline
\end{tabular}
\begin{tabular}{|c|l c|l c|}
\hline
  Element  & $\norm{S_w-S_{w,h}}_{1,h}$  & order & $\norm{p-p_{h}}_{1,h}$  & order  \\
  \hline
  $8^2$  & $1.266 \times 10^{-2}$ &       & $3.211 \times 10^{-1}$ &       \\
  $16^2$ & $6.310 \times 10^{-3}$ & 1.004 & $1.092 \times 10^{-1}$ & 1.556 \\
  $32^2$ & $3.151 \times 10^{-3}$ & 1.002 & $4.503 \times 10^{-2}$ & 1.278 \\
  $64^2$ & $1.575 \times 10^{-3}$ & 1.000 & $2.121 \times 10^{-2}$ & 1.086 \\
  $128^2$& $7.876 \times 10^{-4}$ & 1.000 & $1.046 \times 10^{-2}$ & 1.019 \\
  $256^2$& $3.948 \times 10^{-4}$ & 0.996 & $5.216 \times 10^{-3}$ & 1.004  \\  \hline
  \end{tabular}
\end{center}
\end{scriptsize} \label{Error_Table1}
\end{table}
\end{example}

\begin{example}\label{example2}
We set the level set function $L(x,y)$ as
\begin{eqnarray*}
L(x,y)=(x-0.5)^2+(y-0.5)^2-\frac{1}{16}.
\end{eqnarray*}
Permeability is  $K^+=1$ and $K^-=0.001$.
The exact solutions are given as:
\begin{align*}
p(x,y,t) =& \left\{
  \begin{array}{ll}
    -10L(x,y)(x-1)(2 - t)K^{-}+100 , & \hbox{$L(x,y) > 0$} \\
    10(1-e^{L(x,y)(x-1)})(2 - t)K^{+} +100 , & \hbox{$L(x,y) \leq 0$,} \\
  \end{array}
\right. \\
S_w(x,y,t) =& \cos x(0.2+0.5t),
\end{align*}
where the source terms $q_{\alpha}$, porosity and viscosities are the same as Example 5.1.
Dirichlet boundary condition on all of $\partial\Omega$ is imposed.

Table \ref{Error_Table2} shows the $L^2$, $H^1$  errors of  saturation, pressure and Darcy velocity when  $T=1$.
The time step $\triangle t$ is chosen as  $ \triangle t=16/n^2 =(64/\pi^2) h^2$.
%

\begin{table} \caption{MFVM based on broken $Q_1$-nonconforming IFEM methods for  {Example} \ref{example2}}
\begin{scriptsize}
\begin{center}
\begin{tabular}{|c|l c|l c|l c|}
  \hline
  Element  & $\norm{S_w-S_{w,h}}_{L_2(\Omega)}$  & order & $\norm{p-p_{h}}_{L_2(\Omega)}$ & order & $\norm{\mathbf{u}-\mathbf{u}_{h}}_{L_2(\Omega)}$  & order  \\
  \hline
  $8^2$  & $3.913 \times 10^{-3}$ &       & $4.192 \times 10^{-2}$ &       & $1.097 \times 10^{-3}$ &       \\
  $16^2$ & $1.145 \times 10^{-3}$ & 1.773 & $9.202 \times 10^{-3}$ & 2.188 & $5.418 \times 10^{-4}$ & 1.018 \\
  $32^2$ & $3.137 \times 10^{-4}$ & 1.867 & $2.155 \times 10^{-3}$ & 2.095 & $2.232 \times 10^{-4}$ & 1.279 \\
  $64^2$ & $8.727 \times 10^{-5}$ & 1.846 & $5.436 \times 10^{-4}$ & 1.987 & $1.082 \times 10^{-4}$ & 1.045 \\
  $128^2$& $2.586 \times 10^{-5}$ & 1.755 & $1.295 \times 10^{-4}$ & 2.069 & $5.162 \times 10^{-5}$ & 1.067 \\
  $256^2$& $8.871 \times 10^{-6}$ & 1.543 & $3.147 \times 10^{-5}$ & 2.041 & $2.523 \times 10^{-5}$ & 1.033 \\
  \hline
  \end{tabular}
\begin{tabular}{|c|l c|l c|}
  \hline
  Element  & $\norm{S_w-S_{w,h}}_{1,h}$  & order & $\norm{p-p_{h}}_{1,h}$ & order  \\
  \hline
  $8^2$  & $4.535 \times 10^{-2}$ &       & $6.446 \times 10^{-1}$ &   \\
  $16^2$ & $2.260 \times 10^{-2}$ & 1.005 & $2.187 \times 10^{-1}$ & 1.559 \\
  $32^2$ & $1.124 \times 10^{-2}$ & 1.008 & $1.104 \times 10^{-1}$ & 0.987 \\
  $64^2$ & $5.619 \times 10^{-3}$ & 1.000 & $5.173 \times 10^{-2}$ & 1.093 \\
  $128^2$& $2.836 \times 10^{-3}$ & 0.986 & $2.561 \times 10^{-2}$ & 1.014 \\
  $256^2$& $1.489 \times 10^{-3}$ & 0.930 & $1.275 \times 10^{-2}$ & 1.006 \\  \hline
  \end{tabular}
\end{center}
\end{scriptsize} \label{Error_Table2}
\end{table}

\end{example}

Example 5.1 and Example 5.2  show second order convergence in $L^2$-norm and first order in $H^1$-norm for the pressure variable,
and first order in $L^2$-norm for the velocity variable.
At least one and half  order is observed in $L^2$-norm and first order in $H^1$-norm for the saturation.

\begin{example}\label{example3}
We set the level set function $L(x,y)$ as
\begin{eqnarray*}
L(x,y)=(2-x-y).
\end{eqnarray*}
The exact solutions are given as:
\begin{align*}
p(x,y,t) =& \left\{
  \begin{array}{ll}
    -10L\cos(l_{xy})(2 - t)K^{-}+100 , & \hbox{$L(x,y) > 0$} \\
    -10L\cos(l_{xy})(2 - t)K^{+} +100 , & \hbox{$L(x,y) \leq 0$,} \\
  \end{array}
\right. \\
S_w(x,y,t) =& \left\{
  \begin{array}{ll}
    (1 - 4(l_{xy}+0.25(l_{xy}-2)^2)K^{-})(1 - 0.5t) , & \hbox{$L(x,y) > 0$} \\
    (1 - 8l_{xy}- 4(l_{xy}+0.25(l_{xy}-2)^2)K^{+})(1 - 0.5t) , & \hbox{$L(x,y) \leq 0$,}
  \end{array}
\right. ,
\end{align*}
where the source terms $q_{\alpha}$, porosity and viscosities are the same as Example 5.1. The entry pressure of capillary pressure is $p_d=1$.
Dirichlet boundary condition on all of $\partial\Omega$ is imposed.

We test this example with two different case of permeability.
Table \ref{Error_Table3} shows the $L^2$, $H^1$  errors of  saturation, pressure and Darcy velocity  when  $T=1$  with permeability ($K^+=0.02$, $K^-=0.008$) and Table \ref{Error_Table4} shows errors when $T=1$ with permeability ($K^+=0.1$, $K^-=0.001$).
The time step $\triangle t$ is chosen as  $ \triangle t=16/n^2 =(64/\pi^2) h^2$.

The case of permeability ($K^+=0.02$, $K^-=0.008$) shows second order in $L^2$-norm and first order in $H^1$-norm for the pressure variable, and first order in $L^2$-norm for the velocity variable.
At least one and half order is observed in $L^2$-norm and around one half order in $H^1$-norm for the saturation.

The case of permeability ($K^+=0.1$, $K^-=0.001$) shows similar convergence rates for the pressure and velocity variables.
The average convergence rates for the saturation is $1.527$ in $L^2$-norm  and $0.455$  in  the $H^1$-norm.
Here, an average convergence rate of $L^2$-norm is computed by
$$\frac{log(\norm{S_w-S_{w,h_3}}_{L_2(\Omega)}/\norm{S_w-S_{w,h_8}}_{L_2(\Omega)})}{log(h_3/h_8)} , $$
where $h_{k}= \pi/2^{k+1}$. An average convergence rate of $H^1$-norm is computed similarly.
\end{example}

\begin{table} \caption{MFVM based on broken $Q_1$-nonconforming IFEM methods for  {Example} \ref{example3} with permeability ($K^+=0.02$, $K^-=0.008$)}
\begin{scriptsize}
\begin{center}
\begin{tabular}{|c|l c|l c|l c|}
  \hline
  Element  & $\norm{S_w-S_{w,h}}_{L_2(\Omega)}$  & order & $\norm{p-p_{h}}_{L_2(\Omega)}$ & order & $\norm{\mathbf{u}-\mathbf{u}_{h}}_{L_2(\Omega)}$  & order \\
  \hline
  $8^2$  & $6.569 \times 10^{-4}$ &       & $8.301 \times 10^{-3}$ &       & $1.293 \times 10^{-4}$ &       \\
  $16^2$ & $2.178 \times 10^{-4}$ & 1.592 & $3.581 \times 10^{-3}$ & 1.213 & $6.142 \times 10^{-5}$ & 1.074 \\
  $32^2$ & $6.543 \times 10^{-5}$ & 1.735 & $9.465 \times 10^{-4}$ & 1.920 & $2.415 \times 10^{-5}$ & 1.346 \\
  $64^2$ & $2.141 \times 10^{-5}$ & 1.611 & $2.390 \times 10^{-4}$ & 1.985 & $1.079 \times 10^{-5}$ & 1.163 \\
  $128^2$& $6.022 \times 10^{-6}$ & 1.830 & $5.985 \times 10^{-5}$ & 1.998 & $5.208 \times 10^{-6}$ & 1.050 \\
  $256^2$& $2.476 \times 10^{-6}$ & 1.282 & $1.485 \times 10^{-5}$ & 2.011 & $2.580 \times 10^{-6}$ & 1.013 \\   \hline
\end{tabular}
\begin{tabular}{|c|l c|l c|}
  \hline
   Element  & $\norm{S_w-S_{w,h}}_{1,h}$  & order & $\norm{p-p_{h}}_{1,h}$ & order  \\
  \hline
  $8^2$  & $7.526 \times 10^{-3}$ &       & $4.059 \times 10^{-2}$ &       \\
  $16^2$ & $5.300 \times 10^{-3}$ & 0.505 & $1.857 \times 10^{-2}$ & 1.128 \\
  $32^2$ & $3.567 \times 10^{-3}$ & 0.571 & $7.968 \times 10^{-3}$ & 1.221 \\
  $64^2$ & $2.682 \times 10^{-3}$ & 0.412 & $3.803 \times 10^{-3}$ & 1.067 \\
  $128^2$& $1.842 \times 10^{-3}$ & 0.542 & $1.874 \times 10^{-3}$ & 1.021 \\
  $256^2$& $1.304 \times 10^{-3}$ & 0.498 & $9.343 \times 10^{-4}$ & 1.004 \\ \hline
  \end{tabular}
\end{center}
\end{scriptsize} \label{Error_Table3}
\end{table}

\begin{table} \caption{MFVM based on broken $Q_1$-nonconforming IFEM methods for  {Example} \ref{example3} with permeability ($K^+=0.1$, $K^-=0.001$)}
\begin{scriptsize}
\begin{center}
\begin{tabular}{|c|l c|l c|l c|}
  \hline
  Element  & $\norm{S_w-S_{w,h}}_{L_2(\Omega)}$  & order & $\norm{p-p_{h}}_{L_2(\Omega)}$ & order & $\norm{\mathbf{u}-\mathbf{u}_{h}}_{L_2(\Omega)}$  & order \\
  \hline
  $8^2$  & $2.636 \times 10^{-3}$ &       & $1.738 \times 10^{-2}$ &       & $1.069 \times 10^{-4}$ &       \\
  $16^2$ & $1.282 \times 10^{-3}$ & 1.040 & $3.972 \times 10^{-3}$ & 2.130 & $4.442 \times 10^{-5}$ & 1.267 \\
  $32^2$ & $3.132 \times 10^{-4}$ & 2.033 & $9.503 \times 10^{-4}$ & 2.063 & $1.687 \times 10^{-5}$ & 1.396 \\
  $64^2$ & $2.260 \times 10^{-4}$ & 0.471 & $2.272 \times 10^{-4}$ & 2.064 & $7.385 \times 10^{-6}$ & 1.192 \\
  $128^2$& $3.328 \times 10^{-5}$ & 2.763 & $5.819 \times 10^{-5}$ & 1.965 & $3.535 \times 10^{-6}$ & 1.063 \\
  $256^2$& $1.326 \times 10^{-5}$ & 1.328 & $1.434 \times 10^{-5}$ & 2.021 & $1.748 \times 10^{-6}$ & 1.016 \\   \hline
\end{tabular}
\begin{tabular}{|c|l c|l c|}
  \hline
  Element  & $\norm{S_w-S_{w,h}}_{1,h}$  & order & $\norm{p-p_{h}}_{1,h}$ & order  \\
  \hline
  $8^2$  & $5.553 \times 10^{-2}$ &       & $1.619 \times 10^{-1}$ &       \\
  $16^2$ & $4.702 \times 10^{-2}$ & 0.240 & $6.134 \times 10^{-2}$ & 1.400 \\
  $32^2$ & $2.669 \times 10^{-2}$ & 0.817 & $2.824 \times 10^{-2}$ & 1.119 \\
  $64^2$ & $3.201 \times 10^{-2}$ &-0.262 & $1.396 \times 10^{-2}$ & 1.016 \\
  $128^2$& $1.618 \times 10^{-2}$ & 0.984 & $6.991 \times 10^{-3}$ & 0.998 \\
  $256^2$& $1.147 \times 10^{-2}$ & 0.497 & $3.501 \times 10^{-3}$ & 0.998 \\ \hline
  \end{tabular}
\end{center}
\end{scriptsize} \label{Error_Table4}
\end{table}

\begin{example}\label{example4}
In this example, we use the settings of five spot flooding problem in \cite{bastian1999numerical} (chapter 7.2.) for domain, parameters $\Phi, \lambda_{\alpha}$ and initial/boundary conditions,
but we use different heterogeneities in the domain.

We consider parameters
\begin{align*}
\left\{
\begin{array}{ll}
\Phi=0.2, & \\
\mu_w=0.001[Pa \cdot s] & \mu_n=0.02 [Pa \cdot s] \\
\rho_w=1000[kg/m^3] & \rho_n=1000[kg/m^3] \\
K^+=10^{-10}[m^2], & K^-=10^{-14}[m^2], \\
L(x,y)= ((x-85)^2+(y-185)^2 )^{0.5}-50, &
\end{array}
\right.
\end{align*}
with initial condition,
\begin{align*}
S_w(x,y,0)=
\left\{
\begin{array}{ll}
0.8, &  {\rm if }  (x,y) \in [30,140]\times[170.7,243.3] \\
0, & {\rm otherwise.}
\end{array}
\right.
\end{align*}

Figure 7 shows  $S_w$ at time $t=120[day]$, $t=240[day]$ and $t=375[day]$.
The time step  is chosen as $ \triangle t=(1/240)h^2[day]$.
%
%

\end{example}

\section{Conclusion}
We have introduced a new IMPES type of numerical method {for} two phase flows through heterogeneous porous media  {based on frameworks of IFEM and MFVM}.
{Our scheme can be implemented to solve problems with heterogeneities of various shapes since we can use grids not necessarily aligning  with the interface.}
The pressure variable is approximated by IFEM using  {Rannacher-Turek} nonconforming element
and Darcy velocity is obtained by the local recovery techniques of MFVM of Kwak \cite{chou2003mixed, kwak2010analysis}.
This is similar to, but different from the mixed hybrid method commonly used in porous media.
We believe the MFVM technique  {in the context of \cite{chou2003mixed}} and  IFEM technique are applied to the porous media problem for first time in this paper.
Saturation equation is solved by the control volume method together with an upwinding scheme.
The numerical results  yield almost optimal orders of convergence in all variables for the zero capillary cases.
For the nonzero capillary problem, convergence rates for the saturation deteriorate  but are still optimal for other variables.

\end{document}